\documentclass[12pt,reqno]{amsart}

\usepackage{comment}

\usepackage[final,
color
]{showkeys}
\definecolor{refkey}{gray}{.85}
\definecolor{labelkey}{gray}{.85}

\usepackage{AKstyle}

\usepackage[pagebackref=true, colorlinks]{hyperref}

\hypersetup{pdffitwindow=true,linkcolor=blue,citecolor=blue,urlcolor=blue,menucolor=blue}

\begin{document}%%

\author{Alex Kontorovich}
\thanks{Kontorovich is partially supported by NSF grant DMS-1802119 and BSF grant 2020119.}
\email{alex.kontorovich@rutgers.edu}
\address{Department of Mathematics, Rutgers University, New Brunswick, NJ}

\title{NOTES ON A PATH TO AI ASSISTANCE IN MATHEMATICAL
REASONING}

\maketitle

These informal notes are based on the author's lecture at the National Academies of Science, Engineering, and Mathematics workshop on ``AI to Assist Mathematical Reasoning'' in June 2023. 
The goal is to think through a path by which we might arrive at AI that is useful for the research mathematician.

\section{Finding a Holy Grail}

I find utility in trying to work backwards: What could be our end goal in this game, our ``Holy Grail''?
One possibility is the following.

\begin{def*}
     AI  solves the Riemann hypothesis.
\end{def*}
Of course this applies more generally to any major, longstanding problem of interest to mathematicians.
 I see two ways this could go:

\begin{enumerate}
    \item AI might give a million-line, dense, incomprehensible proof of RH. In this nightmare scenario, people like me will spend the rest of their lives just trying to understand what it's saying and why.\footnote{Perhaps (1) is anyway unreasonable, as some \textit{other} AI will digest the proof  and explain it to us in terms we understand?}
    \item Alternatively, AI might give a perfectly comprehensible, beautiful proof of RH! (Is this a dream? Or also a nightmare?! Now I'm \textit{really} out of business, and spending my life prompting GPT instead of relishing the thought of solving a hard problem.)
\end{enumerate}

Another potential target, instead of solving \textit{our} problems, may be the following.

\begin{def*}
    AI  starts making its own beautiful definitions, conjectures, and theorems, doing ``Alien math'' way ahead of what humans can comprehend and utilize.
\end{def*}

The history of mathematics is full of vignettes, such as Monstrous Moonshine, or the Dyson-Montgomery tea conversation, in which patterns across seemingly unrelated fields were serendipitously discovered and exploited to uncover (perhaps still conjectural) underlying structures. 
Humans already do this on their own, though such discoveries are relatively few and far between. AI has also already proved somewhat effective here; see, e.g., \cite{WilliamsonLackenby}, and the more recent discovery \cite{HeEtAl} of murmurations in the Fourier coefficients of elliptic curves. (That said, in these examples, ``AI'' really refers to: statistical inference analysis on large datasets produced via human-written deterministic algorithms. In this sense, it is not all that different from Gauss conjecturing the Prime Number Theorem from large tables of primes.) 
%such as the following.  At some point people noticed that the list of Fourier coefficients of the $j$-invariant of elliptic curves seems eerily related to that of the dimensions of irreducible representations of the Monster (the largest sporadic simple group). This is exactly the kind of drunk-on-Moonshine observation I looked for as EiC of the journal Experimental Mathematics. Eventually ideas from vertex operator algebras got involved, and Borcherds got a Fields Medal.
%The point is: it’s rare, but humans \textit{can} do this. Can AI?
Can AI discover structures such as perfectoid spaces or infinity categories? If so, might it find novel  structures, definitions, conjectures, and theorems, with sufficient rapidity that human beings have no hope of keeping up?

As far as I can tell, both of these Holy Grails  are distant science fiction.
Here is a more modest, medium-term possibility.

\begin{def*}
    AI can \emph{assist} mathematicians in theorem proving. 
\end{def*}

Certainly if AI can’t even do that, then there’s no hope of the more ambitious Holy Grails. What might we consider ``success'' here?

The current workflow of a research mathematician typically involves something like:
\begin{itemize}
    \item 
Formulating an idea for a major theorem or a class of theorems.
\item 
Trying to break up these theorems into smaller propositions.
\item Deriving each proposition from a set of even smaller lemmas.
\item At each scale (lemma/proposition/theorem), employing standard techniques or exploring new, non-standard approaches, perhaps doing literature searches or looking for analogous arguments that have appeared in other settings.\footnote{LLMs like ChatGPT are already capable of \textit{sometimes} giving helpful pointers, but this mostly happens  when humans have already found connections and written about them.}
\item
Iterating through the whole process, modifying the statements, refining approaches, and discovering alternative techniques, etc etc, all hoping for eventual success.
\end{itemize}

A lot of this does feel like it could be outsourced. 
In other sciences, PIs have labs of PhD students and postdocs who fight the various local battles in parallel, with the PI serving as a General overseeing the action and making ``global'' adjustments to the plan of attack. For many reasons, we don't/can't do this in math. But could AI help automate some of these steps, perhaps drastically speeding up the workflow of a mathematician capable of harnessing these tools? 

Let's continue working backwards from Holy Grail 3 as our target.

%So what if we had AI help accelerate these parts of the process? What would that entail? 

\section{Necessity of Interactive Theorem Provers}

As experience shows, for ``AI'' to do anything, it first needs to be trained on extremely \textbf{large} data sets. A first question might be: data sets of \textit{what} exactly?

\begin{Conj}\label{conj:1}
    Large language models (LLMs) trained on natural language alone will \emph{not} reason reliably at the level of professional mathematics.
\end{Conj}

Mathematics is communicated in such a way as to  appear to be a language amenable to sampling. But at its core, the underlying arguments are always deterministic: either logically sound or not. Language itself exhibits stochastic properties, with various ways to express ideas and infer meaning. Mathematics, on the other hand, relies on an underlying language, such as English, but its essence is entirely deterministic and precise.

For instance, if an LLM generates a 98\% believable Bach chorale, we would be delighted. However, if it produces a 98\% correct mathematical argument, it may be completely  worthless. That missing 2\% might be a reduction to a problem 
as difficult as the original, or it may just be flat out wrong.
%\\
%
\begin{comment}
    
\textbf{What is ``AI''?}
Perhaps it’s worthwhile for a moment to recall what ``AI'' actually is, as there seems to be quite a lot of confusion in the general public. AI is a fancy word for ``algorithm''. We tell the computer to multiply $2.01$ by $3.69$ and it does it. It doesn’t do anything that we do not explicitly tell it to do. We tell it to do lots and lots of linear algebra and statistics in rather complicated ways, and then we’re surprised by the results. Because the algorithm is so complicated, it’s very hard to predict what the outcome will be.

But it’s \textit{just} an algorithm. Nobody knows why this class of algorithms seems to work so well on so many different kinds of problems. Of course the design was meant to mimic neurons and dendrites, but as far as I know, there is very little rigorous theoretical understanding.
\end{comment}
%
%
%
%
The processes of ``production'' and ``editing'' are distinct functions, even within the human brain. LLMs excel at stochastic generation, producing the next word (or token) by sampling from a statistical distribution,
%\footnote{An aside: Even the ``statistical sampling'' part isn’t entirely true, because getting a computer to do something truly random is itself a difficult problem! So in practice people substitute pseudorandom processes, that is, deterministic processes which seem to behave randomly, often seeded by something ``random'' like the current clock tick.}
 whereas mathematics demands deterministic editing. That is my perspective, for what it's worth.
\\

\textbf{Counterpoint:}
OpenAI, Google, and many other groups much smarter than me and with many more resources than I have are actively working to disprove \conjref{conj:1}! They believe that, if only LLMs were given \textit{enough} training data, enough parameters, enough transformers, etc, 
then they \textit{would} solve professional-level mathematical problems, as they've already succeeded handsomely with a variety of standardized math and reasoning tests. (Their work in this direction is not unrelated to issues of ``alignment'', which we have neither the time nor expertise to delve into here...)

And even if I’m completely wrong, and LLMs are indeed capable of producing, in natural language, something that reads like a perfect math paper, how can we ever trust it? We would be burdened with the responsibility of refereeing the thousands or millions of papers they generate to discern the correct ones from mere ``hallucinations''.

\begin{Conj}
The path to AI assisting research mathematicians is through an adversarial process, likely involving Interactive Theorem Provers (e.g., Lean, Isabelle, Coq, etc).    
\end{Conj}

\section{A refined Holy Grail}

Here is a suggested mechanism to aim for Holy Grail 3. Imagine the following scenario.

\begin{def*}
One asks a ChatGPT-like prompt in natural language about a Lemma or technique idea. The LLM bounces back and forth with Lean,\footnote{Here and throughout, we will use ``Lean'' as a shorthand for any Interactive Theorem Prover, perhaps even one not yet invented.} and eventually outputs in natural language: ``here is what I’m able to prove about your question'', together with a formalized certificate that the argument described is valid.    
\end{def*}

This possibility seems, to me, much less out of reach, at least in the medium-term.
If one accepts this as a desired goal, the question again becomes: how do we train AI for this task? Any of the GPT4, Bard, or other such systems, to work nearly as well as they do,  typically require training on trillions of data points, whether these are measured in bytes, lines of text, tokens, etc. In contrast, if we consider the current formalized libraries of mathematics, we might have on the order of 10 million ``data points'', again, whether that means:  raw bytes, lines of code, pairs of goal states and next proof lines, or pairs of natural language lines and formal lines, or something else entirely.

This naturally leads to the following:

\begin{Conj}
The current rate-limiting step for AI proof assistance  is: producing orders of magnitude more lines of  formalized 
professional-level
mathematics.    
\end{Conj}

If you believe this conjecture, then the next natural question is: how to produce lots more?

\section{Paths to Large Datasets of Formalized Mathematics}

Producing high-quality formalized mathematics is difficult! Modern mathematics requires interactions of many different fields in concert. One needs huge libraries that can all interact with one another compatibly (this is one of the main features of Lean’s ``mathlib'' library).
I see (at least) three reasons to be somewhat optimistic here.

\begin{figure}
    \centering
    \includegraphics[width=4in]{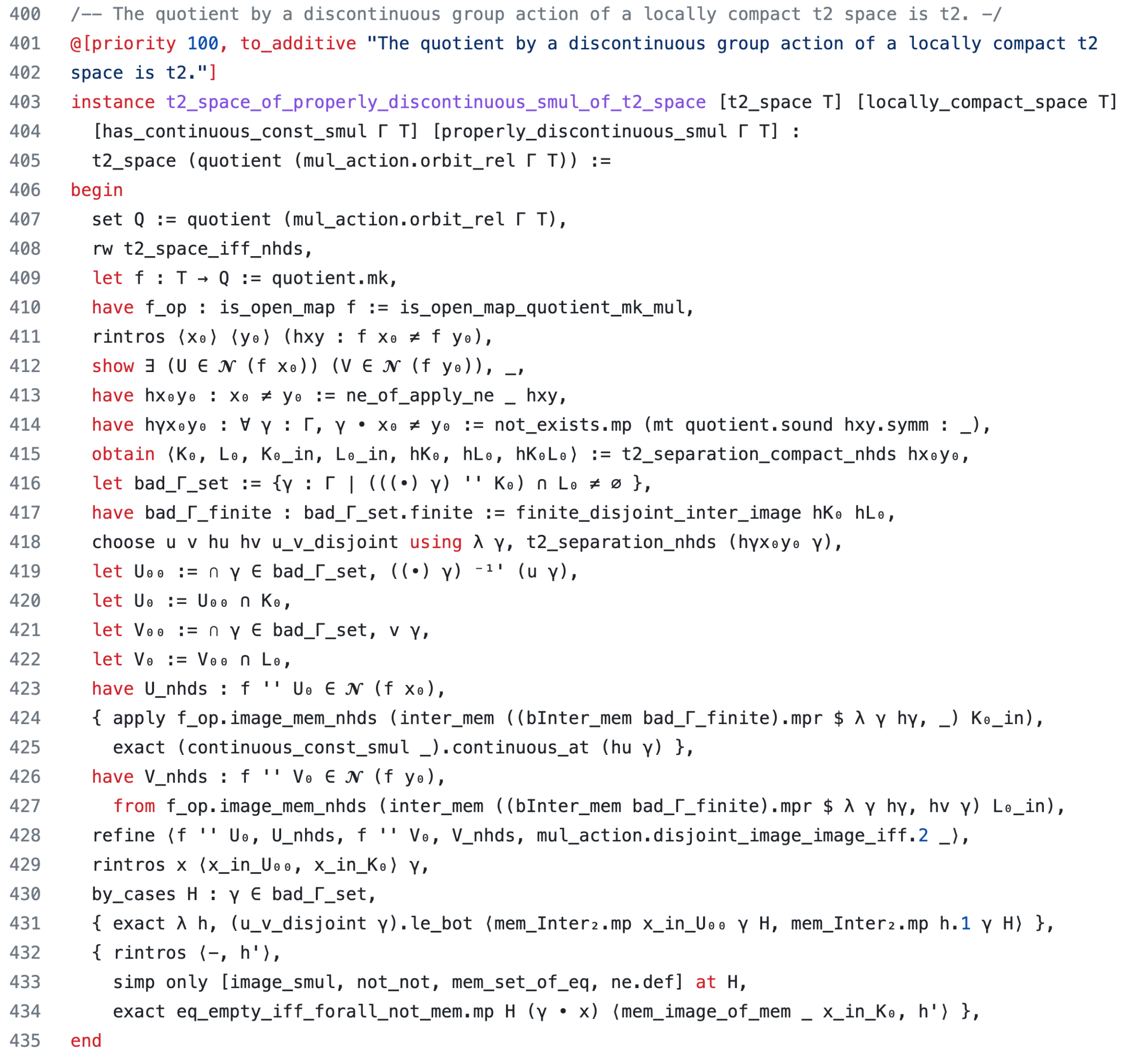}
    \caption{A formal proof from \cite{KontorovichMacbeth2022}}
    \label{fig:1}
\end{figure}

\begin{enumerate}
    \item The pace of production has been steadily increasing. Three vignettes to illustrate this are:
    \begin{enumerate}
        \item It took two years from Ellenberg-Gijswijt's 2017 solution  \cite{EllenbergGijswijt2017} of the Cap Set Conjecture to its 2019 formalization \cite{DahmenHolzlLewis2019}.
        \item It took one year from Gardam's 2021 counterexample \cite{Gardam2021} to Kaplansky's unit conjecture for group rings to its 2022 formalization \cite{GadgilTadipatri2022}.
        \item The formalization \cite{BloomMehta2022} of Bloom's solution \cite{Bloom2021}  to the Edros-Graham density conjecture on Egyptian unit fractions arrived \textit{before} its referee reports!
    
    \end{enumerate}
    (Of course these results are cherry-picked, and it is still relatively difficult to formalize most research-level mathematics...)

    \item The levels of complexity that Theorem Proving software can handle is advancing. We've gone from:
    \begin{enumerate}
        \item Proving difficult statements about simple objects, such the formalization \cite{Avigad2004} of the Prime Number Theorem: the primes are elementary to define, but the statement on their asymptotic behavior is relatively non-trivial to establish.
        \item Proving elementary statements about complicated objects, such as the formalization \cite{BuzzardCommelinMassot2020} of the mere \textit{definition} of perfectoid spaces.
        \item Proving difficult statements about complicated objects, such as the success of the Liquid Tensor Experiment \cite{Scholze2022, LTE2022}.
    \end{enumerate}
    
    \item With enough progress on Application Programming Interface (API), formalized proofs can (sometimes) be made to work almost \textit{exactly} how  human ones do. One example of this is the formalization \cite{KontorovichMacbeth2022} of the author and Heather Macbeth of the fact that the quotient $X/\Gamma$ of a locally compact Hausdorff space $X$ by a discontinuous group action $\G$ is itself Hausdorff. This statement is certainly nothing earthshattering; it's the kind of thing a 1st year PhD student (or advanced undergrad) should be able to solve, but it's also not entirely trivial. What's more remarkable, to me at least, is that the proof, shown in \figref{fig:1}, reads almost verbatim like a human one, with nearly each line of the formalization having a corresponding line in a natural language argument. (Of course one needs to first become familiar with the syntax of Lean, in the same way that one eventually becomes accustomed to reading dollar signs and backslashes in LaTeX and visualizing the outcome of compilation...) 
\end{enumerate}

Here is a big reason to be pessimistic.

\begin{theorem}
All that said, humans alone will \emph{never} reach a trillion lines of formalized mathematics.    
\end{theorem}

The proof is obvious. Equally obvious, then, is that humans need automated assistance in formalization! This idea is far from original, and was proposed already by Szegedy in 2020 \cite{Szegedy2020}. 
Before AI can assist mathematicians  with solving \textit{new} problems, it had better get really good at formalizing known (to humans) solutions!
There are many groups already working hard on this from a wide variety of viewpoints, see, e.g., 
\cite{RabeLeeBansalSzegedy2020, SJLRSJS2022, LLLMHERL2022, JWZLLJLWL2022}
for but a sample. Many more are needed, in this author's opinion.\footnote{Added in print: the authors of \cite{GZACD2023} are able to achieve impressive results in the restricted realm of coding, training on ``only'' billions of tokens, rather than trillions. That said, the path to the former may not be so different from the latter; that is, a mechanism to reach a billion, say, lines of formalized mathematics, seems likely to also reach a trillion.}
\\

In closing, it seems prudent to direct resources towards developing a positive feedback loop between human and AI formalization,
with progress in each driving the other.
This can hopefully eventually build up enough of a training database from which other AI systems can  effectively learn, and become truly useful and reliable as research assistants. Who knows where things will go from there.
\\

\textbf{Acknowledgements:} The author would like to thank Kevin Buzzard, Drew Sutherland, and Geordie Williamson for many comments and suggestions that improved on an earlier draft.

\bibliographystyle{alpha}

\bibliography{AKbibliog}

\end{document}